\documentclass[12pt,a4paper]{report} 
\usepackage{latexsym,amsfonts,theorem,rotating,amsmath} 

\theoremstyle{break}
\theorembodyfont{\slshape}
\newtheorem{Thm}{Theorem}[section]

\newtheorem{Lem}[Thm]{Lemma}
\newtheorem{Cor}[Thm]{Corollary}

\theorembodyfont{\upshape}
\newtheorem{Def}[Thm]{Definition}
\newtheorem{Rem}[Thm]{Remark}

\newtheorem{Exs}[Thm]{Examples}
\newtheorem{Question}[Thm]{Question}

\newcommand{\BoP}[1]{\noindent {\sc Proof#1: }} 
\newcommand{\EoP}{\hfill$\Box$\vspace{6pt}}	    
\newcommand{\EoPP}{\EoP\vspace{0.5cm}}          
\newcommand{\bb}{\mathbb}

\newcommand{\N}{{\bb N}}

\begin{document}
\pagestyle{plain}
\pagenumbering{arabic}
\linespread{1.03}

\setcounter{chapter}{3}

\chapter{Deformation of cylinder knots}
\section{Introduction}
In Chapters 1 - 3 we studied Lissajous knots and their diagrams,
symmetric unions and cylinder knots.\footnote{These chapters
correspond to articles \cite{Lam2}, \cite{Fourier}, \cite{L3} and \cite{LO}.} 
We learned that a cylinder knot with coprime parameters $n$ and $m$ is a symmetric union.

In this chapter we define a family $R(n_1,n_2,n_3)$ of Lissajous
knots which are symmetric unions. Then we succeed in establishing
a relationship between Lissajous knots and cylinder knots:
if certain conditions are satisfied the partial knot of $R(n_1,n_2,n_3)$
is the same as the partial knot of $Z(2n_1,n_2,n_3)$.
Although these knots are different in general, by Theorem 2.6 in \cite{L3}
their determinants are equal.

The second topic of this chapter is the possibility to exchange
the parameters $n$ and $m$ of cylinder knots.
The examples
$Z(3,11,13)=Z(3,13,11)$ and 
$Z(3,11,16)=Z(3,16,11)$, but
$Z(4,11,13)\not = Z(4,13,11)$
lead to the

\begin{Question}\label{question1}
Under which circumstances is $Z(s,n,m)$ equal to $Z(s,m,n)$?
\end{Question}

To answer this question and to analyze the first problem we consider
billiard knots in a solid torus.
This solid torus we choose as a cube with identified front and back face
(a ``flat solid torus'').
The billiard knots in it have the parameters $s,n,m$ as known from 
cylinder knots. But here $n$ and $m$ are exchangeable.

We develop a procedure to deform billiard curves in a cylinder to
billiard curves in a flat solid torus. During the deformation
singularities can occur and the knot type can change.
We can express the result of our investigation in a simplified way
like this: if both billiard curves have the property that their
knot type is the same at the beginning and the end of the deformation
(we say the curve is ``weakly stable'') then $Z(s,n,m)=Z(s,m,n)$.

\section{Lissajous knots which are symmetric unions}
We start with the definition of the knots $R(n_1,n_2,n_3)$ which are
a subfamily of Lissajous knots (or billiard knots in a cube).
As in Chapter 3 let $g$ be the sawtooth function 
$g(t)=2|t- \lfloor t \rfloor -\frac{1}{2}|$.

\begin{Def}
Let $n_2$ be odd and $R(n_1,n_2,n_3)$ be the billiard knot in
a cube defined by the map
$t \mapsto (g(n_1t+\frac{1}{4}),
g(n_2t),g(n_3t+\gamma)),\,\,t\in[0,1]$, with $\gamma > 0$ sufficiently
small.
\end{Def}

The meaning of the condition ``sufficiently small'' will become
clear in the following lemma:
the phase $\gamma$ is moved away from zero so that 
the singularities in the knot are resolved and no new
ones are added.

\begin{Lem}
The knots $R(n_1,n_2,n_3)$ are symmetric unions.
\end{Lem}

\BoP{}
Because of $g(\frac{1}{4})=\frac{1}{2}$ and $g(0)=1$, the projection
on the $x$-$y$-plane is symmetric to the axis $\{(x,y)|x=\frac{1}{2}\}$.
If $\gamma=0$ there is a maximum at $(\frac{1}{2},1)$ and the crossings
on the symmetry axis are singular.
Lemma 1.2 in \cite{Jones} shows that there are no other
singularities.
A small deformation from $\gamma=0$ to $\gamma > 0$ yields a symmetric
union.
\EoP

\bigskip
\begin{figure}[htbp]
\centerline{\includegraphics[height=5cm]{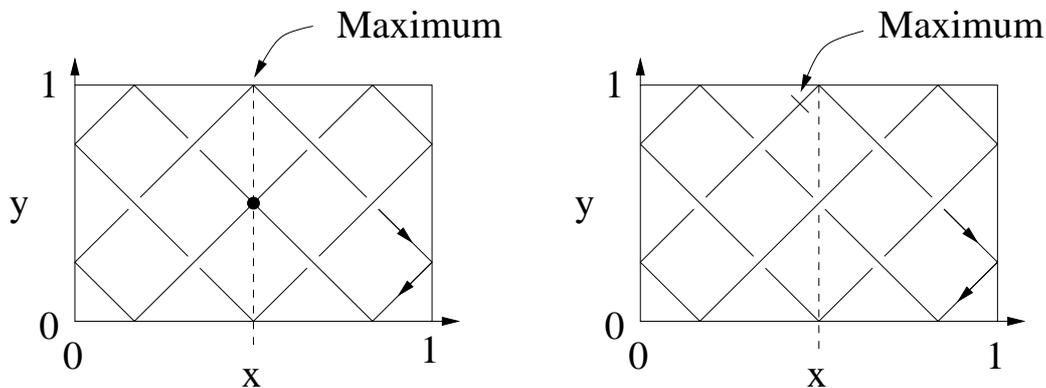}}
\caption{\small The construction of billiard knots in a cube
which are symmetric unions.} 
\label{billsym}
\end{figure}

For some examples of knots from the families $Z(2s,n,m)$ and $R(s,n,m)$
we computed the determinants. For $s=2$ and $n=11$ we get the following
results:

\bigskip
\begin{center}
\begin{tabular}{c|c|c|}
$m$ & $\det Z(4,11,m)$ & $\det R(2,11,m)$ \\
\hline
35 & 1        & $99^2$ \\
37 & 1        & $571^2$\\
39 & $113^2$  & $113^2$\\
41 & 1        & $329^2$\\
43 & $187^2$  & $187^2$\\
\hline
\end{tabular}
\end{center}

\bigskip
As we mentioned in the introduction\footnote{Not part of this chapter.} these numbers made us ask the 

\begin{Question}\label{question2}
Is there a relationship between the families $Z$ and $R$? 
\end{Question}

The two question \ref{question1} and \ref{question2} started the
analysis in this chapter and they were the motivation 
to prove the determinant formula in Theorem 2.6 in \cite{L3}.
Our answers are contained in the Corollaries \ref{switch} and \ref{ZR}.

\section{Billiard knots in a flat solid torus}
\begin{Def}
We call $VT={\bf I}^3/(0,y,z)\sim(1,y,z)$, the cube with identified front and
back face, a \emph{flat solid torus} and
the periodic billiard curves in it  
\emph{billiard knots in a flat solid torus}.
\end{Def}

A notation for billiard knots in a flat solid torus is 
\begin{eqnarray*}
T(s,n,m,\phi):[0,1]&\rightarrow & VT\\ 
t &\mapsto & (s \cdot t \bmod 1, g(nt), g(mt+\phi)). 
\end{eqnarray*}

\bigskip
\begin{figure}[htbp]
\centerline{\includegraphics[height=4.5cm]{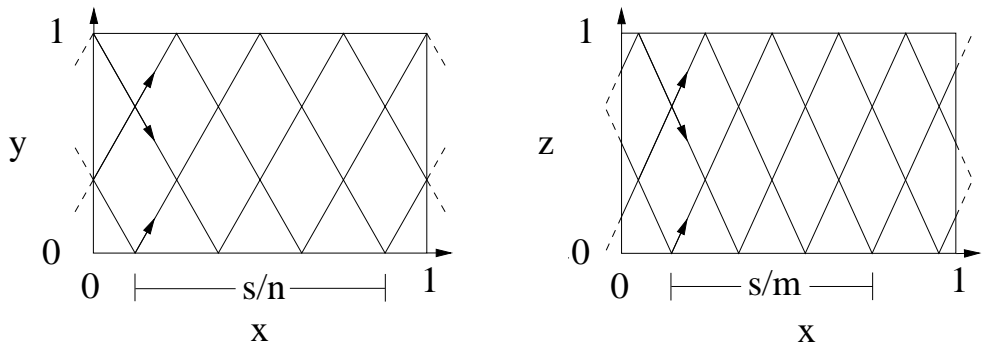}}
\caption{\small The parameters of billiard knot in a flat solid torus.
The figures show the projections on the $x$-$y$-plane and the
$x$-$z$-plane, respectively.} 
\label{}
\end{figure}

The integers $s, n, m$ are the number of strings, the number of 
reflections to the respective faces and $\phi \in [0,1]$ is a phase. 
We need $\gcd(s,n)=\gcd(s,m)=1$ (the parameters
$n$ and $m$ are symmetric to each other) but the 
condition $n \ge 2s+1$ of the cylinder theory is not necessary here. 

\begin{Lem}
(i) There is always a phase $\phi$ such that $T(s,n,m,\phi)$ is without
self-intersections.\\
(ii) The knots $T(s,n,m,\phi)$ are independent of the phase (up to taking
mirror image).
\end{Lem}

\BoP{}
(i): A computation of the crossing parameters shows that if all 
$\phi \in [0,1]$ lead to singular knots we have $\frac{mk}{s} \in {\bf Z}$
with a $k\in \{1,\ldots,s-1\}$. Hence if $\gcd(s,m)=1$ and $s\ge 2$ 
this is impossible. If $s=1$ there are no crossings at all in the projection.
For the proof of (ii) we use (i) and the same symmetry argument as
in Lemma 2.3 in \cite{LO}.
\EoPP

By the first part of the lemma we can neglect the phase
and just write $T(s,n,m)$.
The knots in this chapter are unique only up to mirror image. Therefore
we write $K_1=K_2$, if $K_1$ is equivalent to $K_2$, $-K_2$, $K^*_2$ or 
$-K^*_2$.

The proofs of the next theorems are analogous to 
the respective proofs in the last chapter because all of them
are based on the dihedral symmetry of the projection curve. 
We first treat the aperiodic case:

\begin{Thm}
For pairwise coprime parameters $s, n, m$ the knot $T(s,n,m)$ 
is a ribbon knot and we have $T(s,n,m)=T(s,m,n)$. 
If $n$ or $m$ are even then $T(s,n,m)$ is strongly positive amphicheiral. 
\EoP
\end{Thm}

We return to the general case with $\gcd(n,m)=d\ge 1$.

\begin{Thm} \label{period}
If $\gcd(n,m)=d$ then $T(s,n,m)$ has cyclic period $d$ with linking 
number $s$. 
If $n$ is odd and $m$ is even or
vice versa, $T(s,n,m)$ is strongly positive amphicheiral. The
parameters $n$ and $m$ are exchangeable: we have $T(s,n,m)=T(s,m,n)$.
The factor knot corresponding to the cyclic symmetry 
is $T(s,\frac{n}{d},\frac{m}{d})$. It is a symmetric union.
\EoP
\end{Thm}

\noindent
Further analysis of the situation shows the 

\begin{Thm}\label{cycfree}
If $n$ and $m$ are odd then \\
\hspace*{0.3cm}-if $s$ is even $T(s,n,m)$ has cyclic period 2,\\
\hspace*{0.3cm}-if $s$ is odd $T(s,n,m)$ has free period 2. 
\end{Thm}

\BoP{}
We replace in the parametrization $t$ by $t+\frac{1}{2}$. The result
is $(s(t+\frac{1}{2}) \bmod 1, \frac{1}{2}-g(nt),\frac{1}{2}-g(mt+\alpha))$
because the function $g$ has the property $g(t+\frac{1}{2})=\frac{1}{2}
-g(t)$.

If $s$ is even this shows that the knot is mapped to itself by a rotation
of $\pi$ around the circle $\{(x,y,z)|y=z=\frac{1}{2}\}$. A non-singular
knot is disjoint from this circle, hence it has cyclic period 2.

If $s$ is odd the knot is invariant under the fixpoint-free map 
$(x,y,z)\mapsto(x+\frac{1}{2},\frac{1}{2}-y,\frac{1}{2}-z)$, hence
it is freely periodic with period 2.
\EoP

\begin{Rem}
If $s$ is even, $n$ and $m$ must be odd, so all knots $T(s,n,m)$ with
even $s$ have cyclic period 2.
\end{Rem}

\begin{Rem}
We consider the knots which are parametrized in cylinder coordinates by
\begin{equation}\label{zylko1}
(\varphi(t),r(t),z(t))=(st,3+\cos(nt),\cos(mt+\phi)),\,\,t\in[0,2\pi].
\end{equation}
They are identical to the knots $T(s,n,m,\phi)$:

We deform the function cosine as we did in the case of Lissajous knots
to a piecewise linear function. Hence we can write 
\begin{equation}\label{zylko2}
(\varphi(t),r(t),z(t))=(2\pi \cdot st,1+g(nt),g(mt+\phi)),\,\,t\in[0,1].
\end{equation}
We map the flat solid torus
$KR=\{(\varphi,r,z)|\varphi \in [0,2\pi],1 \le r \le 2, 0 \le z \le 1\}$
by the homeomorphism 
$(\varphi,r,z)\mapsto (\varphi/2\pi,r-1,z)$ to $VT$. 
The image of a knot as (\ref{zylko2}) is the knot in the flat solid
torus $T(s,n,m,\phi)$. An example for (\ref{zylko1}) with $s=3$, $n=7$, $m=5$ 
can be seen in Figure \ref{zylkoabb}.

\begin{figure}[htbp]
\centerline{\includegraphics[height=6cm]{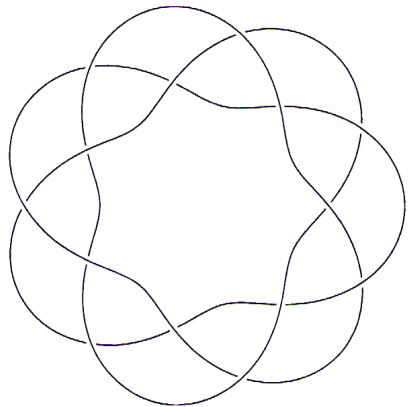}}
\caption{\small The knot $T(3,7,5)$ parametrized in cylinder coordinates.} 
\label{zylkoabb}
\end{figure}
\end{Rem}

\newpage
\section{Stable cylinder knots}
\begin{figure}[htbp]
\centerline{\includegraphics[height=13cm]{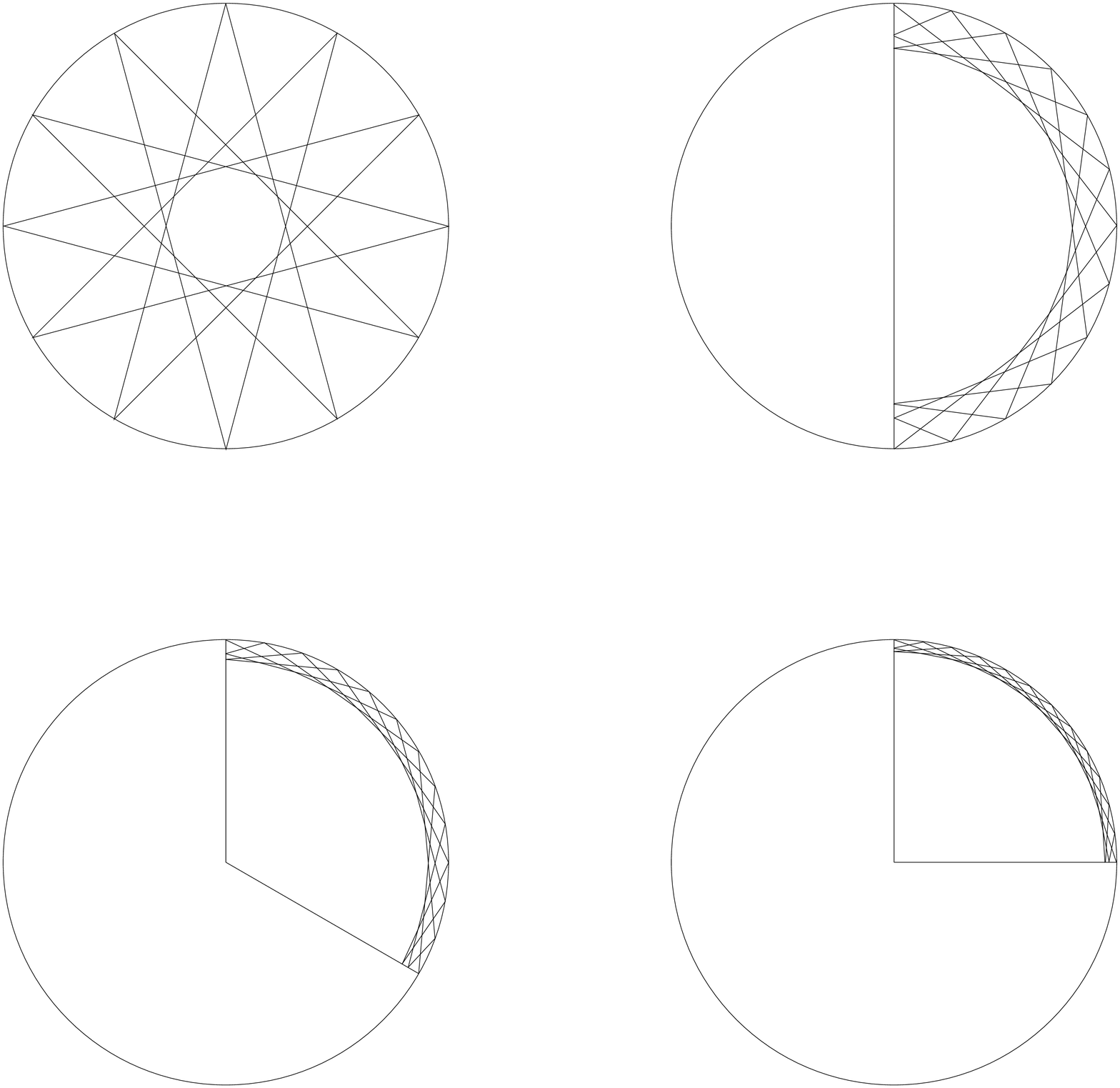}}
\caption{\small The projected billiard curves for $s=5$, $n=12$, $a=1,2,3,4$.} 
\label{def_beta}
\end{figure}

Given an integer $a$ coprime to $s$ 
the knot $Z(s,an,am)$ has cyclic period $a$. The factor knot 
$Z(s,an,am)^{(a)}$ is similar to $Z(s,n,m)$: it consists of 
the same numbers of strings, reflections to the
boundary ${\bf S}^1 \times {\bf I}$ and maxima. 
We expect that these knots ``converge'' in some way for 
$a \rightarrow \infty$ (see Figure \ref{def_beta}). 
Thus if the knots $Z(s,an,am)^{(a)}$ 
are independent of $a$ for large $a$ we write suggestively 
$Z^{st}(s,n,m)=\lim_{a \to \infty} Z(s,an,am)^{(a)}$ for the limit
knot.

Instead of using the factor knot of the periodic knot with period $a$,
we identify the radial faces of a cylinder's $\frac{2\pi}{a}$-slice 
to define $Z(s,an,am)^{(a)}$. 
This description allows us to construct the deformed knots for
$a$ which are not coprime to $s$; we also use the notation $Z(s,an,am)^{(a)}$. 
The condition $n\ge 2s+1$ for cylinder knots is unnecessary for  
the knots $Z^{st}$ because we can choose $a$ so that
$an \ge 2s+1$.

\medskip

\begin{Thm}\label{A}
Let $s,n$ and $m$ be integers with $\gcd(s,n)=1$. Then the knot 
$Z^{st}(s,n,m) := \lim_{a \to \infty}Z(s,an,am)^{(a)}$
is well-defined: there is an integer $\tilde a$  
so that for all $a$ with $a \ge \tilde a$ we have 
$Z(s,\tilde a n,\tilde a m)^{(\tilde a)}=Z(s,an,am)^{(a)}$.
\end{Thm}

We call the knots $Z^{st}(s,n,m)$ \emph{stable cylinder knots}.
We can describe stable cylinder knots in terms of the above defined 
billiard knots in a flat solid torus if the parameters allow it, that
is if $s$ and $m$ are coprime.

\begin{Thm}\label{B}
If the billiard knot in a flat solid torus $T$ is defined for the
parameters $s,n,m$ then $Z^{st}(s,n,m) = T(s,n,m)$. The billiard 
knots in a flat solid torus are a subfamily of the knots $Z^{st}$.
\end{Thm}

\subsubsection*{Answer to Question \ref{question1}}
Now we give an answer to Question \ref{question1}.
A billiard curve in the cylinder with $Z(s,n,m) = Z^{st}(s,n,m)$ 
is called \emph{weakly stable}. In this case we also say that $Z(s,n,m)$
is weakly stable, though this property is not a knot invariant but
a property of the billiard curve.
Our theorems yield the following Corollary \ref{switch}, for
which we assume $\gcd(s,n)=\gcd(s,m)=1$ and $n,m \ge 2s+1$.

\begin{Cor} \label{switch}
If the cylinder knots $Z(s,n,m)$ and $Z(s,m,n)$
are weakly stable then we have $Z(s,n,m) = Z(s,m,n)$. 
\EoP
\end{Cor}

\subsubsection*{Answer to Question \ref{question2}}

Our second Question \ref{question2} is answered by Theorem \ref{C} 
and Corollary \ref{ZR}.

\begin{Thm}\label{C}
If $2s, n, m$ are pairwise coprime, then the knots $T(2s,n,m)$ and
$R(s,n,m)$ are symmetric unions with identical partial knots. Hence
the determinant of $T(2s,n,m)$ equals the determinant of $\det R(s,n,m)$.
\end{Thm}

From this we deduce immediately the

\begin{Cor}\label{ZR}
Let $2s, n, m$ be pairwise coprime and $n \ge 4s+1$. If
$Z(2s,n,m)$ is weakly stable then $\det Z(2s,n,m) = \det R(s,n,m)$. 
\EoP
\end{Cor}

\begin{Rem}
The relationships between the knot families in this chapter
are summarized in the following diagram.

\bigskip
\begin{tabular}{ccccc}
cylinder knots, &&            &&                  \\
$Z(s,n,m)$       &&            && Lissajous knots \\
&&&&\\
$\downarrow$      &&            && $\cup$             \\
&&&&\\
stable cylinder     &&          billiard knots     && Lissajous knots, \\
knots  ,             &$\supset$ & in a flat solid   &$\longleftrightarrow$ & which are symmetric\\
$Z^{st}(s,n,m)$     &&           torus, $T(s,n,m)$ && unions, $R(n_1,n_2,n_3)$ \\
\end{tabular}
\end{Rem}

\bigskip

\subsection{The deformation process}
In the construction of the stable cylinder knots we used slices of
the cylinder with angles $\beta=2\pi/a$. However, to be more 
flexible we now vary the angle $\beta$ continuously 
from $2\pi$ to zero.

We set the total length of the projected curve $\kappa$ onto the 
bottom of the slice to 1 and parametrize $\kappa$ according to its
arclength. The base point $\kappa(0) = \kappa(1)$ is chosen to be a
vertex on the boundary of ${\bf D}^2$. The crossings of $\kappa$ are
numbered from 1 to $(s-1)n$ and the parameters of the $i$-th crossing
are $t_i(\beta)$ and $t'_i(\beta)$.
For the length of $P_0P_b$ in the figure on page 356 in \cite{LO}
depending on $\beta$ we write $\frac{x_{b}(\beta)}{2n}$. Then
$$
t_i(\beta)=\frac{k_i}{2n}+\frac{x_{l_i}(\beta)}{2n},\enspace 
t'_i(\beta)=\frac{k'_i}{2n}-\frac{x_{l_i}(\beta)}{2n},
$$
with $k_i,k_i',l_i \in \N$.
It is remarkable that the sum $t_i(\beta)+t'_i(\beta)$ is 
independently of $\beta$ a constant multiple of $1/2n$.

Let $\Delta_i^\phi(\beta)$ be the height difference of the
billiard curve at the crossing $i$ as a function of $\beta$ and phase
 $\phi\in[0,1]$. Then for $\Delta_i^\phi(\beta)$ we obtain the expression
$$
\Delta^\phi_i(\beta):=g(mt_i(\beta)+\phi)-g(mt'_i(\beta)+\phi).
$$

\subsubsection*{The sign of $\Delta^\phi_i(\beta)$}
Since we are interested
only in the sign of $\Delta^\phi_i(\beta)$ we compare the function $g(t)$
with $h(t)=(\cos(2\pi t)+1)/2$ in the following way in order
to use a theorem from trigonometry:
\begin{eqnarray*}
g(t_1) > g(t_2) &\Leftrightarrow& h(t_1) > h(t_2), \\
g(t_1) = g(t_2) &\Leftrightarrow& h(t_1) = h(t_2), \\
g(t_1) < g(t_2) &\Leftrightarrow& h(t_1) < h(t_2).
\end{eqnarray*}
The formula $\cos(t_1)-\cos(t_2)=-2\sin((t_1+t_2)/2)\sin((t_1-t_2)/2)$ yields
$$
\textrm{sign }\Delta^\phi_i(\beta)=-\textrm{sign }\Bigl(
\sin(\pi[m[t_i(\beta)+t'_i(\beta)]+2\phi])\cdot
\sin(\pi m[t_i(\beta)-t'_i(\beta)])
\Bigr).
$$
This is useful because the first factor is independent of $\beta$ and
the second is independent of $\phi$.

If $\phi$ is chosen so that $\Delta^\phi_i(2\pi)\not=0$, we can form
the quotient
\begin{equation}\label{unabh}
\delta_i(\beta)=\textrm{sign }\frac{\Delta^\phi_i(\beta)}{\Delta^\phi_i(2\pi)}
=\textrm{sign }\frac{\sin(\pi m[t_i(\beta)-t'_i(\beta)])}
                    {\sin(\pi m[t_i(2\pi)-t'_i(2\pi)])}.
\end{equation}
If this sign is $+1$ the sign of the crossing $i$ is the same for the angles
$\beta$ and $2\pi$, if it is $-1$ a crossing change has happened and if it is
zero the crossing is singular for the angle $\beta$.
We stress the fact that the phase $\phi$ is not contained in the expression.

\subsubsection*{Stability of cylinder knots}
After these technical preparations we can define several notions of 
stability.
As in the definition of weakly stable cylinder knots
before Corollary \ref{switch}, the following definition is about
properties of billiard curves in a cylinder, not about knot invariants.

\begin{Def}
The cylinder knot $Z(s,n,m)$ is called
\begin{eqnarray*}
\textsl{weakly stable}                  &:\Leftrightarrow& Z(s,n,m)=Z^{st}(s,n,m),\\
\textsl{positively (negatively) stable} &:\Leftrightarrow& 
      \exists \varepsilon > 0 \textrm{ so that for all } 
      i \textrm{ and for all } \beta \in ]0,\varepsilon]:\\
                               && 
                                 \frac{\Delta^\phi_i(\beta)}
                                 {\Delta^\phi_i(2\pi)}>0 
                                 \,(\,<0\,),\\
\textsl{strongly positive stable}       &:\Leftrightarrow& 
      \textrm{for all $i$ and for all } \beta \in ]0,2\pi]:\\
                               && 
                                 \frac{\Delta^\phi_i(\beta)}
                                 {\Delta^\phi_i(2\pi)}>0. 
\end{eqnarray*}
\end{Def}

If the cylinder knot is positively or negatively stable we call it
also just \emph{stable}. Stable knots are weakly stable and, of course,
strongly positive stable knots are positively stable.

\begin{Exs}
(i) $Z(p,q,q)$ is the torus knot $t(p,q)$ (see \cite{Jones}, p. 153/154). 
It is strongly positive stable because the distribution of $q$ maxima on
the $q$ chords of the projection can be chosen independently of $\beta$.
For instance if we place the maximum and minimum at the ratios $1/4$ and
$3/4$ of the chord's length we get the torus knot
$t(p,q)$ for all $\beta \in ]0,2\pi]$.

(ii) If the parameters $4,n,m$ are pairwise coprime and the knot
$Z(4,n,m)$ is negatively enlaced\footnote{Using the notation of Theorem 4.2 in
\cite{LO} we call the cases $v_1 = -v_3$ and $v_1 = v_3$ {\it negatively enlaced}
and {\it positively enlaced}, respectively.} then it is not stable.
By Theorem \ref{cycfree} we know that the curves $Z^{st}(4,n,m)=T(4,n,m)$ 
are positively enlaced. 
\end{Exs}

\subsection{Proofs of Theorems \ref{A}, \ref{B} and \ref{C}}
\BoP{ of Theorem \ref{A}}
As in Proposition 2.5 in \cite{LO} we compute the parameters $x_l(\beta)$ as
$$
x_l(\beta) = 
\frac{\tan(l\frac{\beta}{2n})}{\tan(s\frac{\beta}{2n})}.
$$
It is not difficult to show that $\frac{d}{d\beta}x_l(\beta)<0$ 
for $\beta \in ]0,2\pi]$, and with de L'Hospital's rule 
we get $\lim_{\beta \to 0}x_l(\beta)=\frac{l}{s}$.

The difference $t_i(\beta)-t'_i(\beta)$ is equal to
$\frac{k_i-k'_i}{2n}+\frac{2x_{l_i}(\beta)}{2n}$ and hence strictly
mono\-tonously decreasing. Hence for all $i$ we find an
$\varepsilon_i > 0$ so that 
sign $\sin(\pi m[t_i(\beta)-t'_i(\beta)])$ is independent of $\beta$
for all $\beta \in ]0,\varepsilon_i]$. For all $a$ with 
$\frac{1}{a}<\min_i \varepsilon_i$ the knots $Z(s,an,am)^{(a)}$ are
the same (up to mirror image, as usual), and we denote them as
$$
\lim_{a \to \infty} Z(s,an,am)^{(a)}=:Z^{st}(s,n,m) 
$$
as we did in the formulation of the theorem.
\EoPP

\begin{figure}[htbp]
\centerline{\includegraphics[height=6cm]{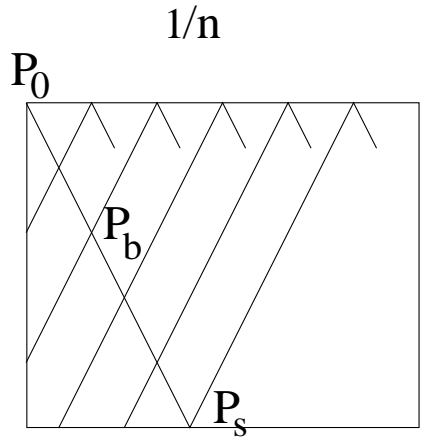}}
\caption{\small The computation of $|P_0P_b|/|P_0P_s|$} 
\label{grenze}
\end{figure}

\BoP{ of Theorem \ref{B}}
The knots $T(s,n,m)$ are exactly the limit case
$\beta = 0$.
To show this we look at Figure \ref{grenze}. There we chose
the same notation as in the figure on page 356 in \cite{LO}.

\bigskip

We have $|P_0P_b|/|P_0P_s|=\frac{b}{s}$. Therefore from
$\lim_{\beta \to 0}x_b(\beta)=\frac{b}{s}$ we deduce that the knot
$T(s,n,m)$ is equal to $Z^{st}(s,n,m)$.
\EoP

\BoP{ of Theorem \ref{C}}
We consider the knots $T(2s,n,m)$ and $R(s,n,m)$ with a maximum pushed to
the axis of symmetry.

\bigskip

\bigskip
\begin{figure}[htbp]
\centerline{\includegraphics[height=5.2cm]{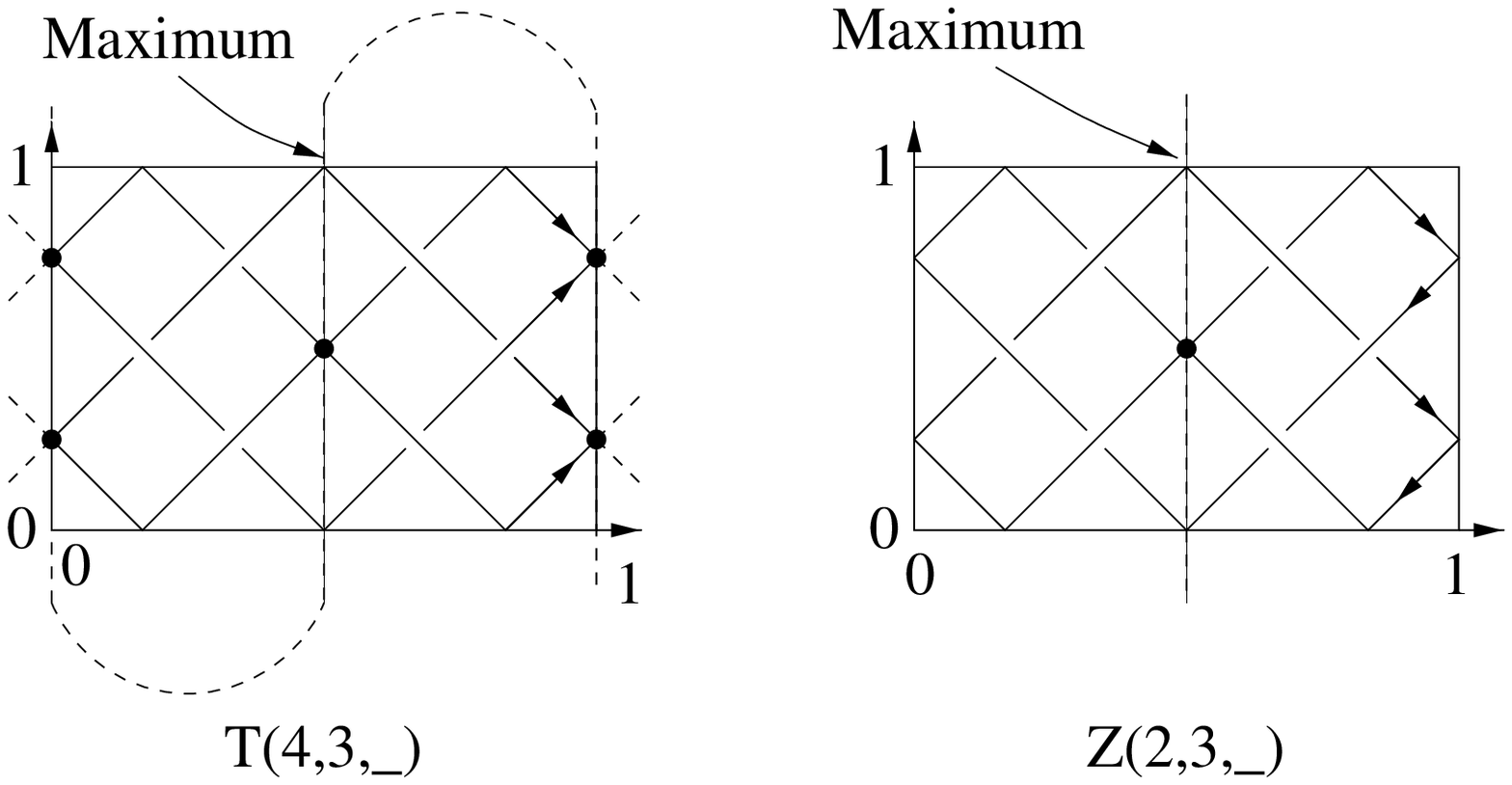}}
\caption{\small Diagrams of the knots $T(2s,n,m)$ and $R(s,n,m)$ 
for $s=2$ and $n=3$ with a maximum on the symmetry axes.} 
\label{beweis3}
\end{figure}

\bigskip
The left and right building blocks of the unions
are the same for $T$ and $R$, because the symmetry axis
of $T$ bends around to the faces of the cube and the reflection in $R$
is exactly the mirrored movement of the straight one in $T$.
\EoP

\subsection{Deformation graphs of cylinder knots}
The difference functions $\Delta^\phi_i(\beta)$ can be displayed graphically
to show whether the cylinder knot is stable or not. If we do this for all
$i$ simultaneously it is better to use (sign $[\Delta^\phi_i(2\pi)])\cdot
\Delta^\phi_i(\beta)$ because then all curves are positive for $\beta=2\pi$
if $\phi$ is non-singular. 
Because of equation (\ref{unabh}) the sign of 
(sign $[\Delta^\phi_i(2\pi)]))\cdot
\Delta^\phi_i(\beta)$ is independent of $\phi$.

We call a graph containing these normalized difference functions
a \emph{deformation graph}.
All we have to do to check the stability of the 
knot is to observe if there is a neighbourhood $]0,\varepsilon]$ of zero
in which all curves are positive or all are negative. Examples are:

\begin{figure}[htbp]
\centerline{\includegraphics[angle=270,width=0.4\textwidth]{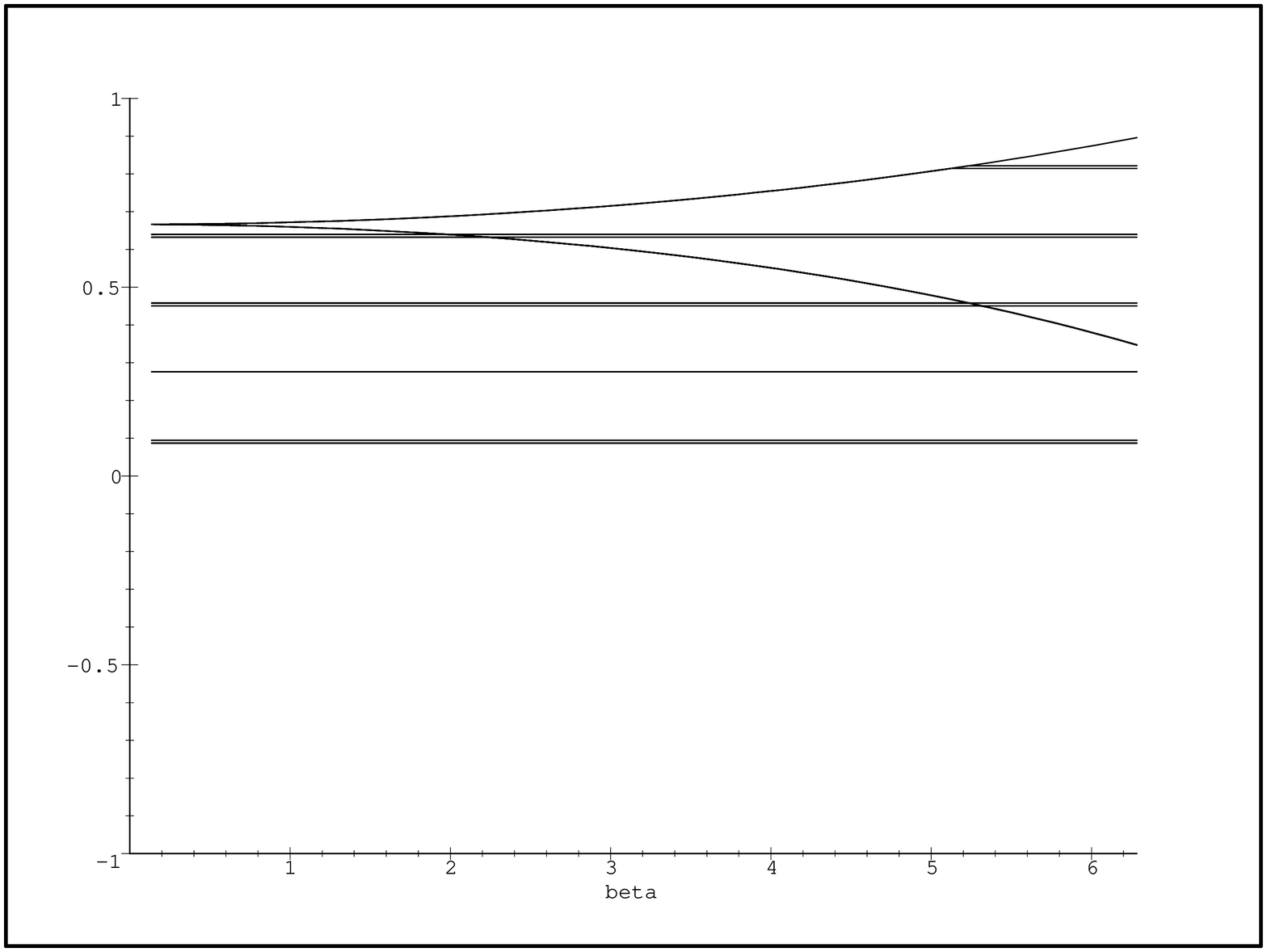} \includegraphics[angle=270,width=0.4\textwidth]{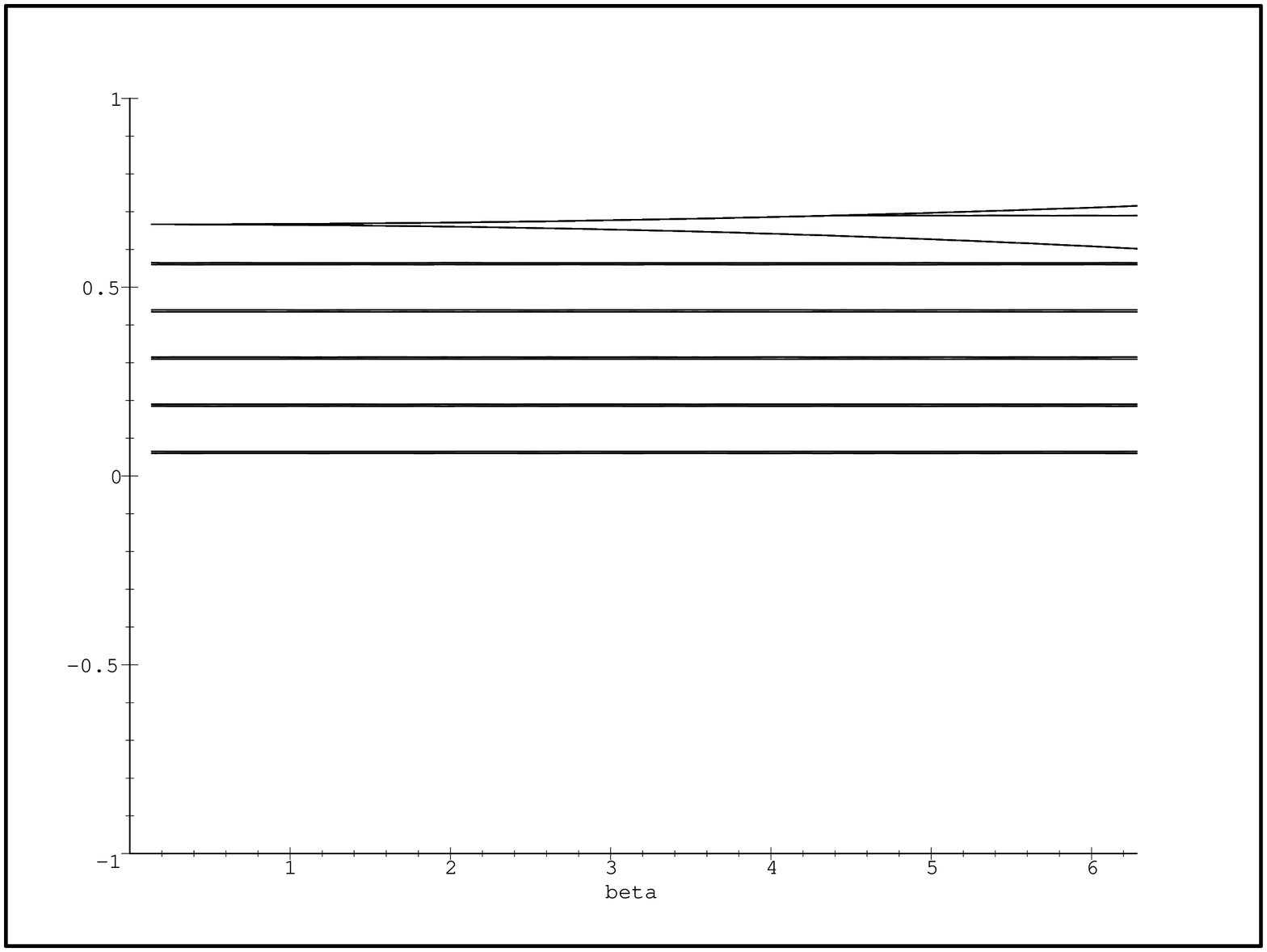}}
\caption{The deformation graphs of $Z(3,11,16)$ and $Z(3,16,11)$} 
\label{dreielf}
\end{figure}

Thus the knots $Z(3,11,16)$ and $Z(3,16,11)$ are
strongly positive stable and Corollary \ref{switch} 
yields $Z(3,11,16) = Z(3,16,11)$.

\begin{figure}[htbp]
\centerline{\includegraphics[angle=270,width=0.4\textwidth]{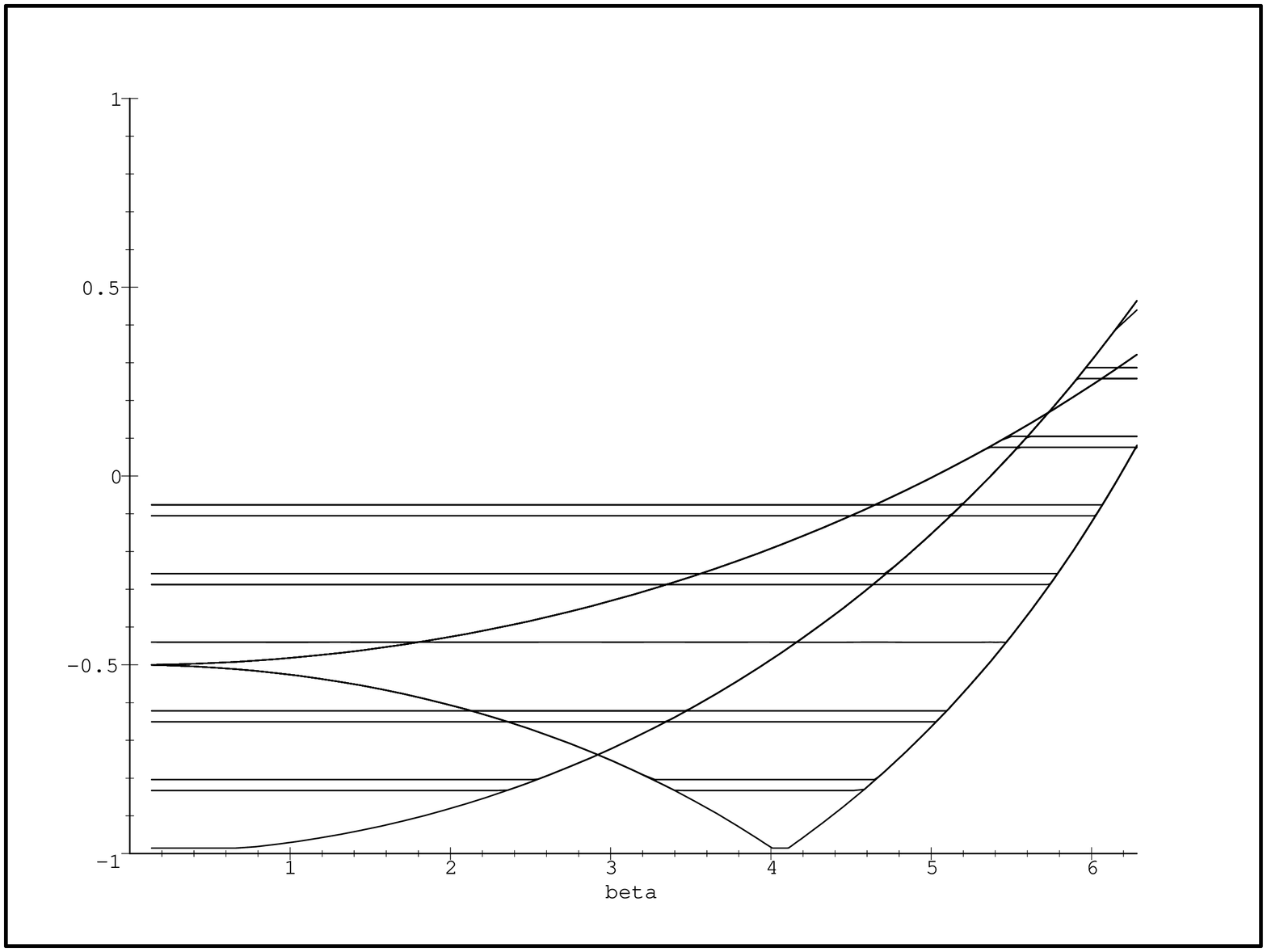} \includegraphics[angle=270,width=0.4\textwidth]{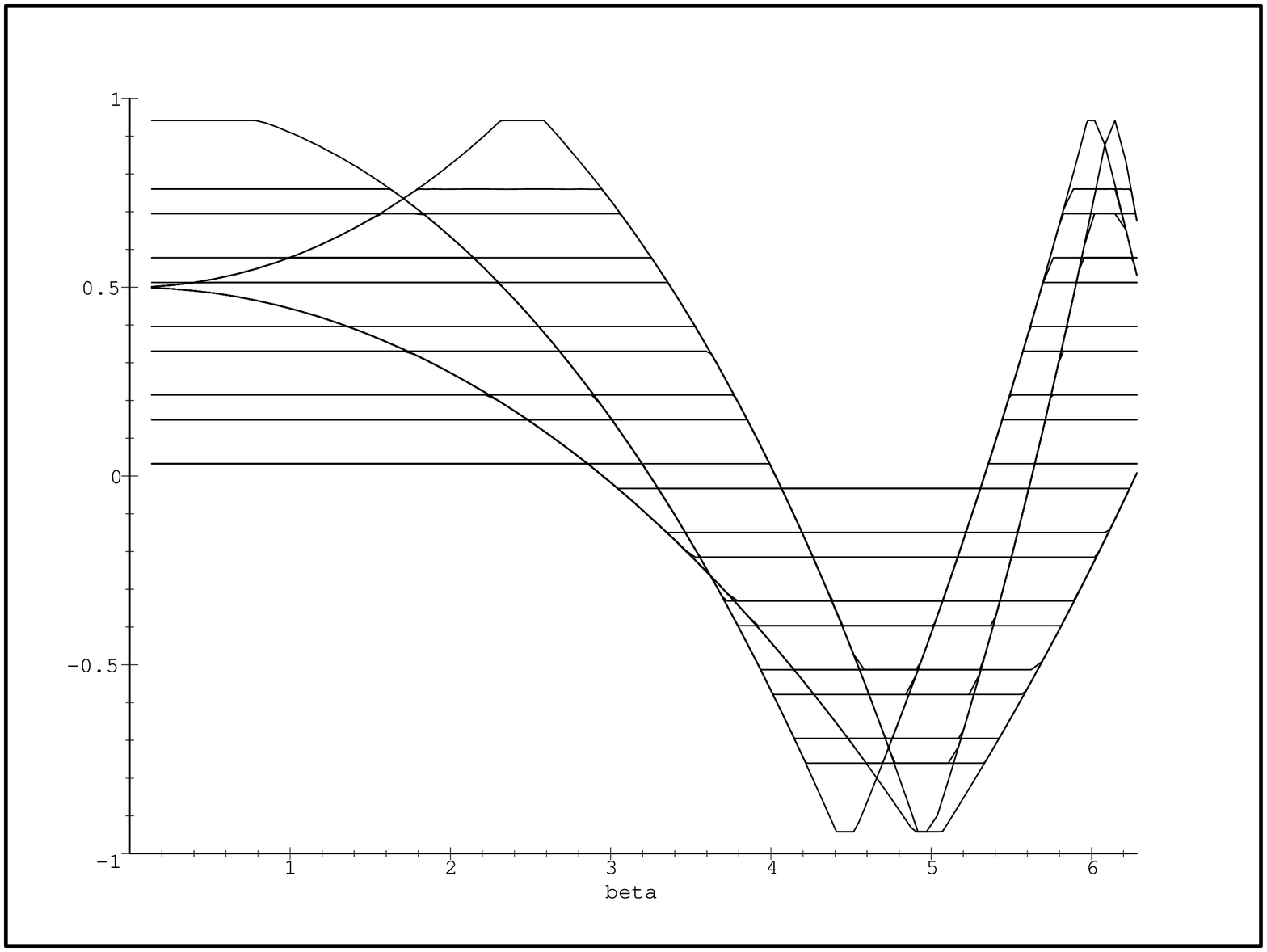}}
\caption{The deformation graphs of $Z(4,11,39)$ and $Z(4,11,119)$} 
\label{vierelf}
\end{figure}

The knot $Z(4,11,39)$ is negatively stable and the knot
$Z(4,11,119)$ is positively stable. Corollary \ref{ZR} 
gives the result $\det Z(4,11,39) = \det R(2,11,39)$.

\begin{figure}[htbp]
\centerline{\includegraphics[angle=270,width=0.4\textwidth]{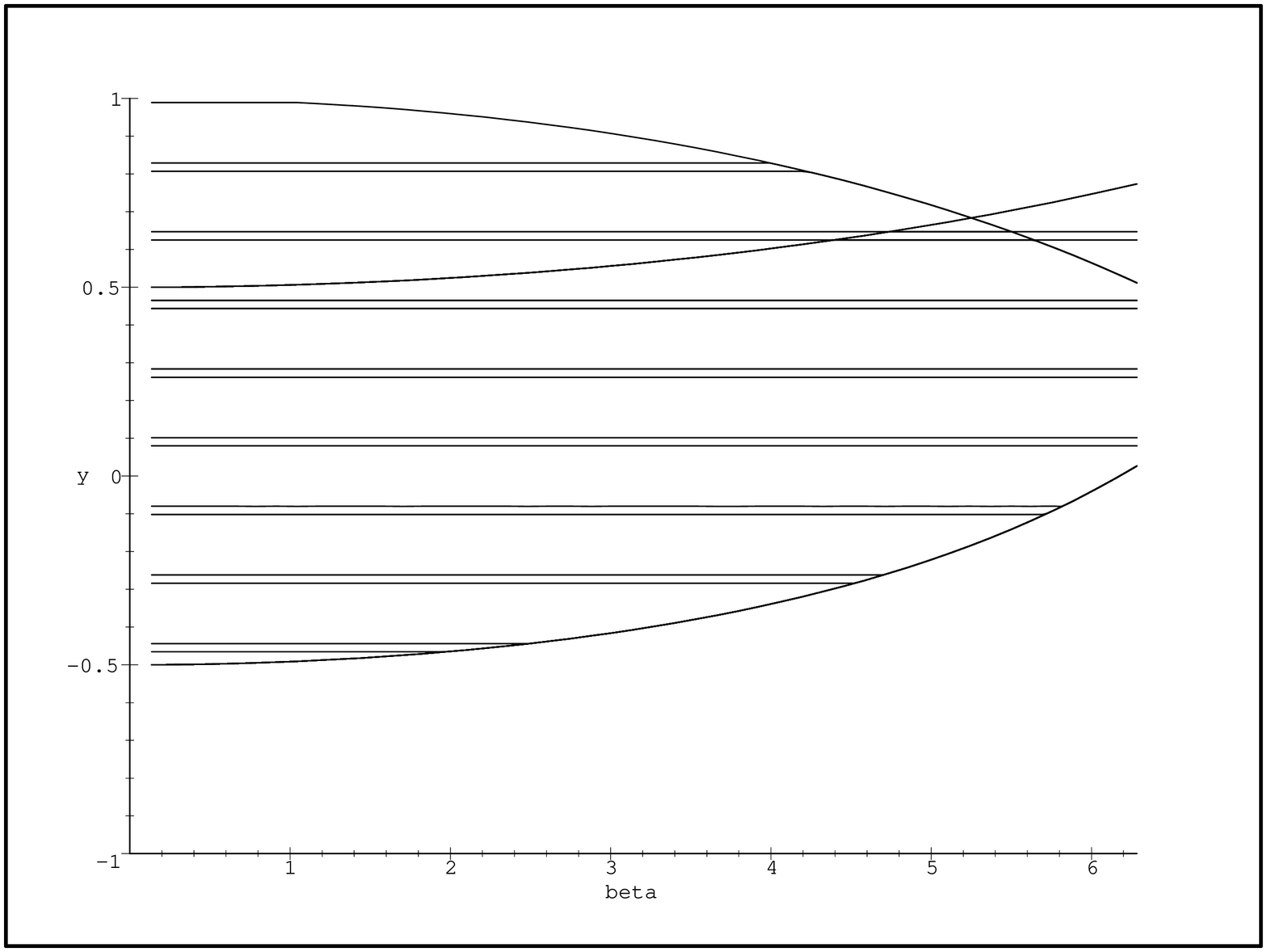} \includegraphics[angle=270,width=0.4\textwidth]{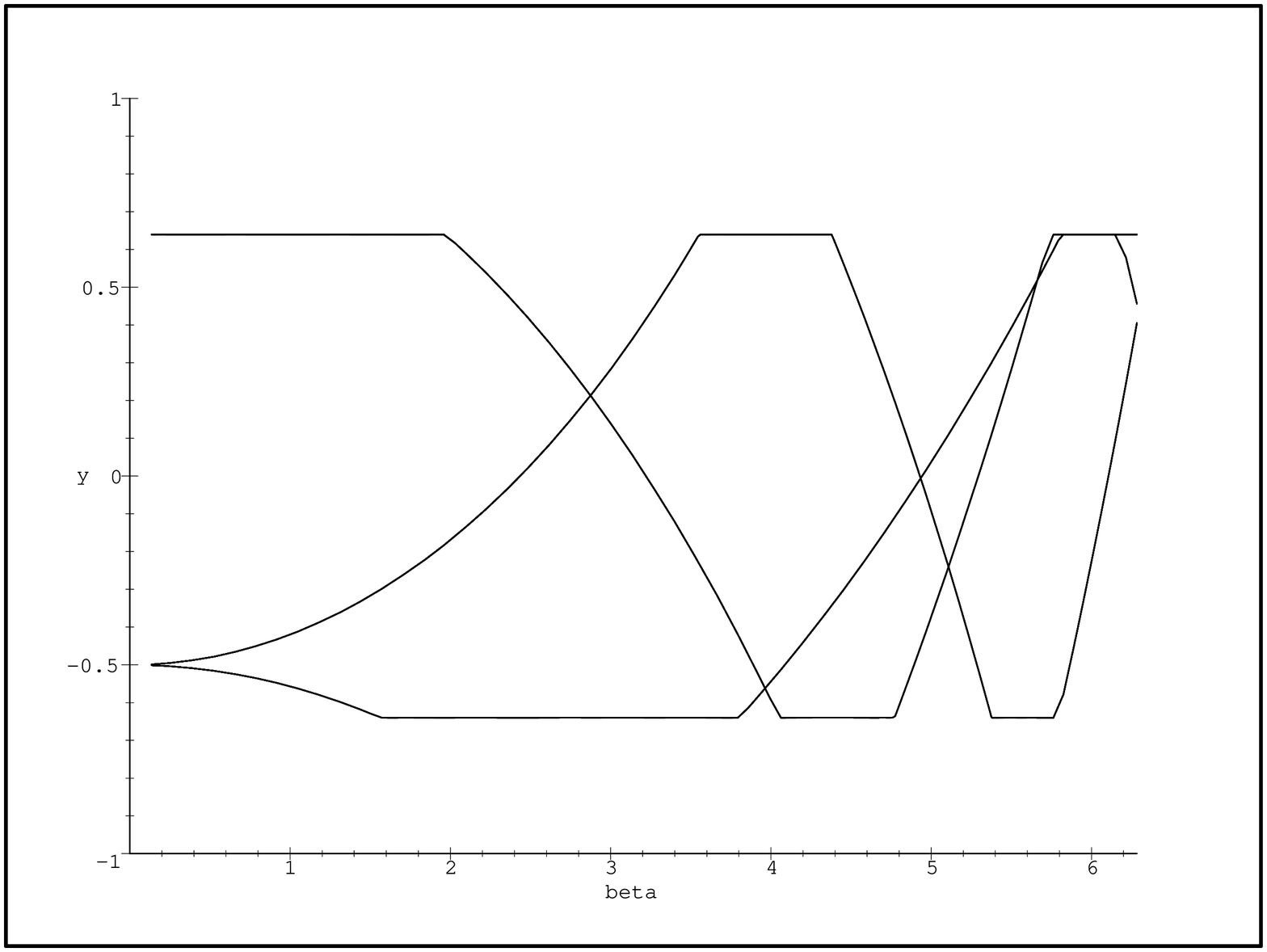}}
\caption{\small The deformation graphs of $Z(4,11,13)$ and $Z(4,11,121)$} 
\label{nichtstabil}
\end{figure}

From Figure \ref{nichtstabil} we get the information that 
$Z(4,11,13)$ and $Z(4,11,121)$ are not stable. The three lines
in the second deformation graph result from the knot's periodicity. 
Therefore only three difference functions occur.

\begin{figure}[htbp]
\centerline{\includegraphics[angle=270,width=0.4\textwidth]{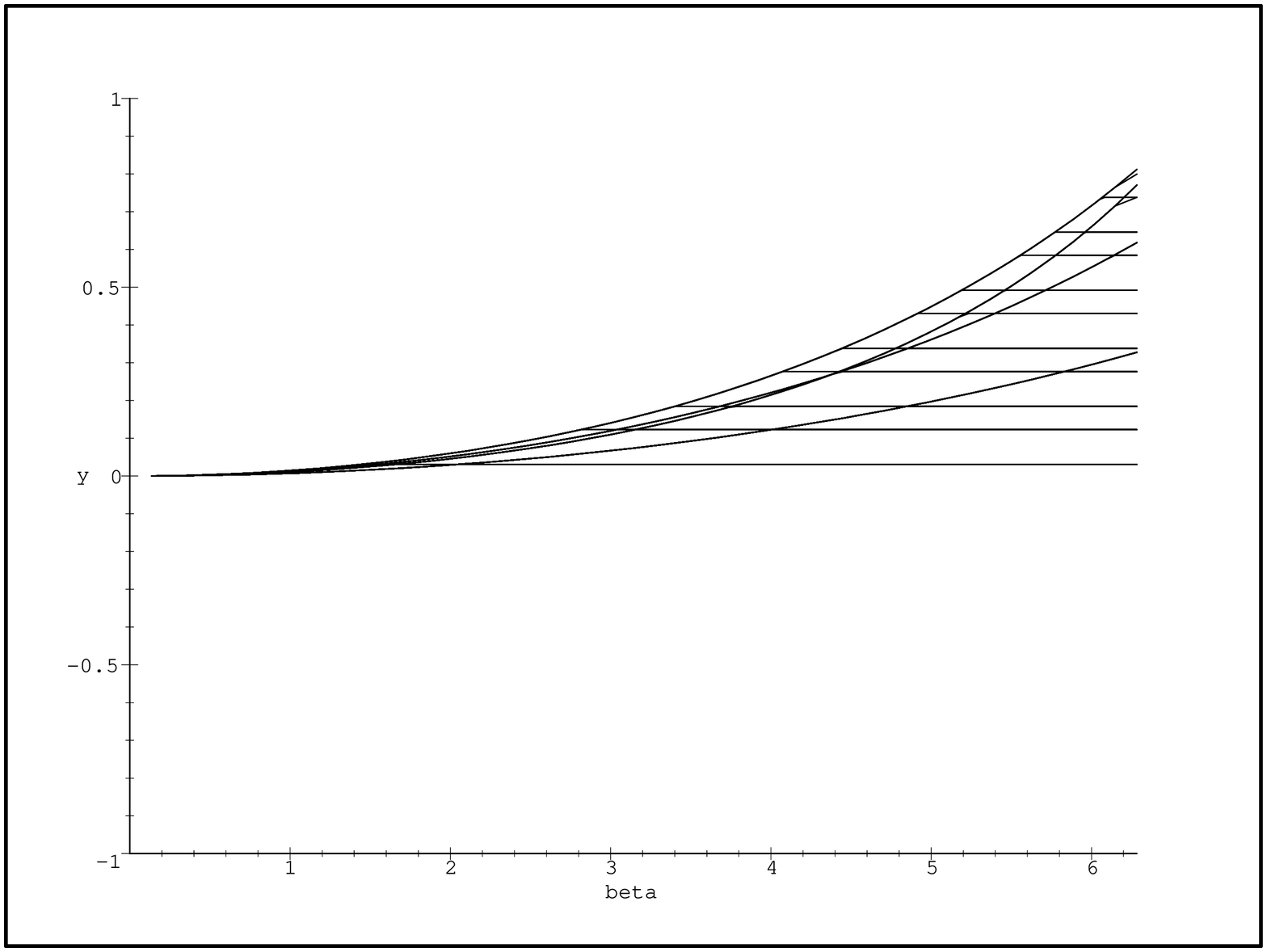} \includegraphics[angle=270,width=0.4\textwidth]{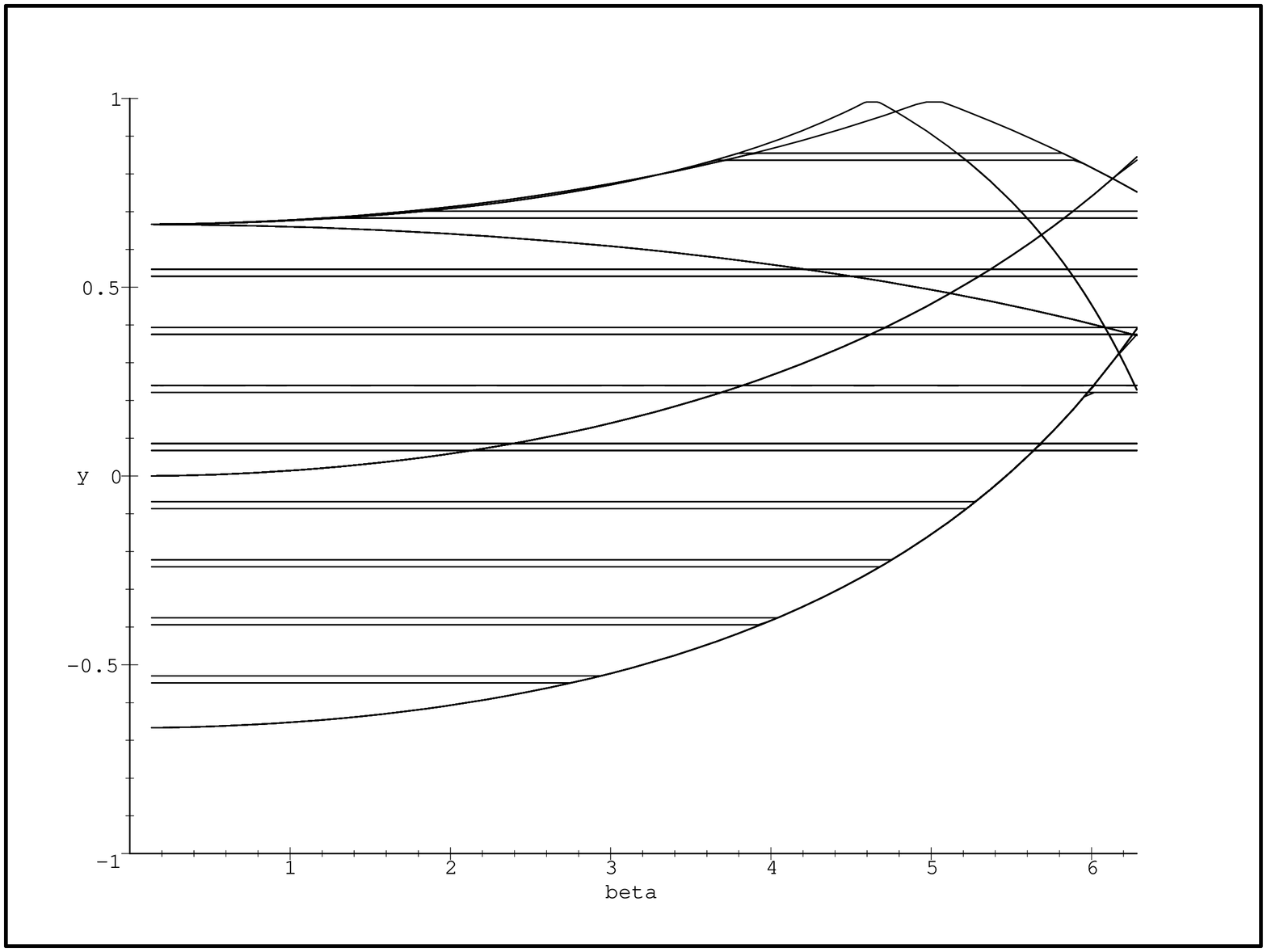}}
\caption{\small The deformation graphs of $Z(5,13,20)$ and $Z(6,13,14)$} 
\label{singu}
\end{figure}

In Figure \ref{singu} examples of deformations are shown
for parameters $s$ and $m$ which are not coprime.
In this case the knot $T(s,n,m)$ does not exist, and there
are difference functions converging to zero for $\beta \to 0$. 
The knot $Z(5,13,20)$ is strongly positive stable and the knot
$Z(6,13,14)$ is not stable.

\bigskip
We give more information on the knots $Z(4,11,m)$ for $m=1,\ldots,156$. 
They are strongly positive stable for
$m=1,\ldots$,12, 14, 15, 16, 18, 19, 20, 22, 23, 24, 28, 32, 36, 40, 44, 48,
positively stable (but not strongly positive stable) for 
$m=$119, 123, 127, 131, 
negatively stable for
$m=$39, 43, 47, 51, 55, 59, 60, 62, 63, 66, 67, 70, 71, 137, 141, 145, 
149, 153
and for the others not stable (they could be weakly stable but we did
not check this.)

\section{Additional parameters of deformation}
\begin{figure}[htbp]
\centerline{\includegraphics[height=8cm]{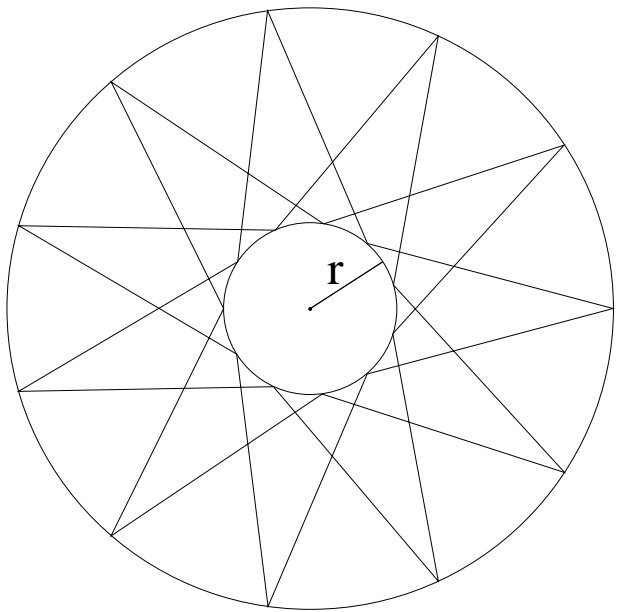}}
\caption{\small The projection of a billiard knot in a prism over
an annulus with $s=3$ and $n=11$.} 
\label{defv}
\end{figure}

In this section we mention two further possibilities to deform
cylinder knots.
For the first one we replace the cylinder by an annulus $\times$ interval.
The inner radius $r$ is variable, see Figure \ref{defv}.

\enlargethispage{2cm}

Another possibility is the study of prisms over ellipses.
The excentricity $\varepsilon$ and the starting point of the
billiard curve are additional parameters. 
By the classical Theorem of Poncelet every point on the ellipse is
starting point of a periodic billiard curve with $s$ rotations and
$n$ reflections, if one point has this property.
Furthermore all such periodic curves have the same length.

\vspace{0.5cm}
\noindent
Added in 2012: 

\medskip
\noindent
Daniel Pecker: {\it Poncelet's theorem and billiard knots}, Geometriae Dedicata {\bf 161}, 323--333 (2012).

\vspace{5cm}
\noindent
Christoph Lamm,\\
R\"{u}ckertstr.~3,\\
65187 Wiesbaden, Germany,\\ 
e-mail: christoph.lamm@web.de

\end{document}